\documentclass[journal]{IEEEtran}
\usepackage{cite}
\usepackage{amsmath,amssymb,amsfonts}
\usepackage{algorithmic}
\usepackage{graphicx}
\usepackage{textcomp}
\usepackage[ruled,vlined]{algorithm2e}
\usepackage[affil-it]{authblk}
\usepackage{blindtext}
\usepackage{CJK}
\usepackage{type1cm}
\usepackage{times}
\usepackage[marginal]{footmisc}
\usepackage[affil-it]{authblk}
\usepackage{color}
\usepackage{multirow}
\usepackage{comment}
\usepackage{bbm}

\usepackage{amssymb}
\usepackage{nomencl}
\usepackage{multirow}
\makenomenclature
\usepackage{etoolbox}
\renewcommand\nomgroup[1]{%
  \item[\bfseries
  \ifstrequal{#1}{A}{Abbreviations}{%
  \ifstrequal{#1}{B}{Variables}{%
  \ifstrequal{#1}{C}{Subscripts and Superscripts}{}}}%
]}

\begin{document}

\title{HVAC Scheduling under Data Uncertainties: A Distributionally Robust Approach}
\author{Guanyu~Tian,~\IEEEmembership{Student Member,~IEEE,}
        Qun~Zhou,~\IEEEmembership{Member,~IEEE,}
        Samy~Faddel,~\IEEEmembership{Member,~IEEE,}
        Wenyi Wang 
        \thanks{\indent This material is based upon work supported by the U.S. Department of Energy's Office of Energy Efficiency and Renewable Energy (EERE) under the Building Technology Office, BENEFIT 2019 Award Number DE-EE0009152. Disclaimer: The views expressed herein do not necessarily represent the views of the U.S. Department of Energy or the United States Government.
        
        \indent Guanyu Tian, Qun Zhou, Samy Faddel, and Wenyi Wang are with the Department of Electrical and Computer Engineering, University of Central Florida, Orlando, FL 32816 USA (e-mail: tiang@knights.ucf.edu; qun.zhou@ucf.edu; samy.faddel@ucf.edu; wenyi.wang@ucf.edu)
        }
        }

\date{\today}

\maketitle
\begin{abstract} 
The heating, ventilation and air condition (HVAC) system consumes the most energy in commercial buildings, consisting over 60\% of total energy usage in the U.S. Flexible HVAC system setpoint scheduling could potentially save building energy costs. This paper first studies deterministic optimization, robust optimization, and stochastic optimization to minimize the daily operation cost with constraints of indoor air temperature comfort and mechanic operating requirement. Considering the uncertainties from ambient temperature, a Wasserstein metric-based distributionally robust optimization (DRO) method is proposed to enhance the robustness of the optimal schedule against the uncertainty of probabilistic prediction errors. The schedule is optimized under the worst-case distribution within an ambiguity set defined by the Wasserstein metric. The proposed DRO method is initially formulated as a two-stage problem and then reformulated into a tractable mixed-integer linear programming (MILP) form. The paper evaluates the feasibility and optimality of the optimized schedules for a real commercial building. The numerical results indicate that the costs of the proposed DRO method are up to 6.6\% lower compared with conventional techniques of optimization under uncertainties. They also provide granular risk-benefit options for decision making in demand response programs. 
\end{abstract}

\begin{IEEEkeywords}
HVAC system, distributionally robust optimization, stochastic optimization, Time of Use rate.
\end{IEEEkeywords}

\IEEEpeerreviewmaketitle

\section{Introduction}
\IEEEPARstart{T}{he} heating, ventilation and air conditioning (HVAC) system, on average, consumes 44\% of the total energy of commercial buildings \cite{eia}. HVAC is the largest contributor of commercial building energy consumption, and its control is incentivized by demand response programs. Through time-varying electricity rates, the economic performance of commercial buildings can be improved while providing benefits to the grid \cite{lu2012evaluation,goddard2014model}. 

This paper studies a flexible HVAC scheduling scheme that plans daily setpoint in the day before the actual operation, aiming to minimize the energy cost for the day ahead. The setpoint scheduling has been applied to the control of residential air conditioning system by smart thermostats and has been proved to be economical \cite{nest,ecobee,honeywell}.

However, any HVAC schedule made in the day before should not compromise occupant comfort in real-time operations. Extra caution is required given that uncertainties exist in HVAC scheduling optimization, including data uncertainty, model uncertainty, and building state uncertainty \cite{tian2018review}. 

This paper delves deep into the impact of data uncertainty, in particular, the impact of ambient temperature on HVAC scheduling, energy cost, and occupant comfort when responding to Time-of-Use (TOU) rate. Many works have been focused on continuously controlling individual components of an HVAC system. The classic method is the proportional-integral-derivative (PID) control that minimizes the error signal through feedback loop using PID controllers \cite{moore2003pump,mirinejad2008control,naidu2011advanced}. Nevertheless, they fail to minimize the total energy consumption of the entire HVAC system in a coordinated manner. Optimal control methods are proposed to fill this gap \cite{engdahl2004optimal,wang2008supervisory,levenhagen1993hvac}. Though the energy efficiency of such optimal controls are higher than the PID methods, in practice, they face robustness issues that are caused by measurement uncertainties. Stochastic optimization methods are then proposed to solve this problem by explicitly taking the uncertainties into account \cite{qiu2019stochastic}. Furthermore, the moment-based distributionally robust optimization methods are incorporated to provide robust solutions \cite{du2017energy}. 

In this paper, the HVAC scheduling problem is first formulated in a deterministic manner. The goal is to minimize the daily energy cost while satisfying the constraints of indoor air temperature comfort and mechanical operating requirement. Then considering the forecasting uncertainty in the ambient temperature, the problem is formulated into a Stochastic Programming (SP) problem and a Robust Optimization (RO) problem. To overcome limitations in SP and RO formulations, and to further enhance the robustness of the optimized schedule against prediction error uncertainty, a Wasserstein metric-based Distributionally Robust Optimization (DRO) approach is then proposed. The ambiguity set is a Wasserstein ball with the predicted probability distribution being the center with a pre-defined radius. The proposed DRO method is initially formulated as a two-stage problem and then reformulated into a tractable mixed-integer linear programming (MILP) form. The DRO-based approach is then applied to a real commercial building and the results show improvement of cost savings and operational robustness. 

The contributions of this paper are three folds:  
\begin{enumerate}
\item \textbf{Incorporating the probability distribution uncertainty in the HVAC scheduling problem through DRO formulation with the  Wasserstein metric.} The proposed method enhances the robustness of the optimal solution by finding the worst probability distribution within an ambiguity set. The worst distribution considered in the Wasserstein metric-based ambiguity set is not based on any assumptions of underlying probability distributions, such as the mean value, symmetricity, and skewness, and hence, the result is expected to be more robust. 
\item \textbf{Reformulating the proposed DRO formulation into a tractable MILP problem.} The initial formulation of the proposed DRO is an intractable two-stage optimization problem, where the fist-stage minimizes the total energy cost and the second-stage finds the worst probability distribution. The tractable form is achieved by applying duality theorem and equivalent substitution of minimization problems with inequality constraints.
\item \textbf{Validating the feasibility, optimality and robustness of DRO derivation using data from a real commercial building.} The one-day schedule performance of the proposed DRO method is compared to DO, SP, and RO methods under 1000 scenarios. Numerical results show that the costs of the proposed DRO method are up to 6.6\% lower compared with conventional techniques of optimization under uncertainties. The DRO approach also provides granular risk-benefit options for decision making in demand response programs. 
\end{enumerate}

This paper is organized as follows. Section \ref{do} introduces the deterministic formulation of the commercial building HVAC scheduling problem as the baseline model. The deterministic formulation is then extended to stochastic optimization including SP and RO formulations by considering the uncertainty of outdoor air temperature in section \ref{soro}. The proposed DRO method and its tractable reformulation are introduced in section \ref{dro}. Section \ref{case} presents the numerical results of the two case studies. Section \ref{conclusions} provides conclusions and discusses future research directions.

\section{Deterministic Formulation of HVAC Scheduling}\label{do}
The HVAC scheduling problem is initially formulated into a Deterministic Optimization (DO) problem using mixed-integer linear programming (MILP) as the decision variables are binary schedules. The overall HVAC scheduling problem is formulated as follows: 
\begin{subequations}
\begin{align}
&\min_{x} \sum_{t=1}^{T}c_t\Delta tP_{t,h}\label{doobj}\\
s.t. \quad &x_{t+1}\geq x_t - x_{t-1},\forall t\geq2\label{updown1}\\
&x_{t+2}\geq x_t - x_{t-1},\forall t\geq2\label{updown2}\\
&x_{t+3}\geq x_t - x_{t-1},\forall t\geq2\label{updown3}\\
&1 - x_{t+1}\geq x_{t-1} - x_t,\forall t\leq T-1\label{updown4}\\
&1 - x_{t+2}\geq x_{t-1} - x_t,\forall t\leq T-2\label{updown5}\\
&1 - x_{t+3}\geq x_{t-1} - x_t,\forall t\leq T-3\label{updown6}\\
&\text{Cooling: }T^{in}_t\leq T^{ub}_t,\forall t \label{summer}\\
&\text{Heating: }T^{in}_t\geq T^{lb}_t,\forall t \label{winter}\\
&T^{in}_t = b_1x_t +b_2T^{oa}_t + b_{3}T^{in}_{t-1} + b_0,\forall t \label{zonetemp}\\
&P_t = a_{1}x_{t} + a_{2}T^{oa}_t + a_0,\forall t \label{powerconsumption}
\end{align}
\label{doformulation}
\end{subequations}
The objective function \eqref{doobj} is to minimize the total operation cost of an HVAC system by optimally scheduling the on/off modes $x$ of the building, where $c_t$ and $P_t$ indicate the electricity rate and power consumption of the building at the $t^{th}$ time interval. $c_t$ is the TOU rate provided by utility companies. $P_t$ is modeled in \eqref{powerconsumption} using linear regression, where $b_1$, $b_2$ and $b_3$ and $b_0$ are the coefficients obtained by training building operational data. $\Delta t$ denotes the duration of a time interval and $T$ is the total number of time intervals. The two status of commercial building HVAC systems are the occupied (on) mode and unoccupied (off) mode, which are represented by the decision variable $x$. 

Frequent changes of setpoint would cause mechanic issues and reduce the life span of HVAC components. Hence, the optimized HVAC schedules should meet the minimum up/down time requirement. \eqref{updown1}-\eqref{updown6} formulates this constraint assuming the sampling interval is 15 minutes and the minimum up/down time is 1 hour \cite{tian2020optimal}.

The functionality of HVAC systems is to maintain the indoor air temperature within the comfort zone, which is a fixed or time-variant range of temperature. The constraints on the indoor air temperature $T^{in}_t$ is formulated by \eqref{summer} and \eqref{winter}, where $T^{ub}_t$ and $T^{lb}_t$ denotes the upper and lower bound of the comfort zone at the $t^{th}$ time interval. \eqref{summer} is applied in summer when HVAC systems operate under cooling conditions and \eqref{winter} is applied in winter when HVAC systems are under heating conditions. Without loss of generality, in this paper, we focus on cooling conditions. The formulations and solutions can be easily extended to heating conditions. The model of indoor air temperature $T^{in}$ is the autoregressive model with exogenous inputs (ARX) formulated in \eqref{zonetemp}, where $a_{1}$ and $a_{2}$ are the coefficients and $a_0$ is the constant term of intercept \cite{hao2016transactive}.

\section{Considering the Uncertainty of Ambient Temperature}\label{soro}
The deterministic formulation lacks robustness because the schedules are computed based on the deterministic prediction of day-ahead outdoor air temperature. Once the actual temperature is different from the prediction, the feasibility and optimally of the schedule are compromised. To improve the robustness of optimal schedules, the outdoor air temperature is treated as a stochastic variable in stochastic optimization formulations. Below we study the most common formulations using Stochastic Programming (SP) and Robust Optimization (RO).

\subsection{Stochastic Programming}\label{so}
In SP, the probability distributions of stochastic variables are used to generate scenarios, and an optimal solution is found across all scenarios. There are variations in SP formulations, yet the scenario-based approach is straightforward and useful \cite{geng2019data}. 

In our application, each scenario is a one-day temperature profile. The objective function of SP is to minimize the weighted average energy cost of all scenarios as formulated in \eqref{spobj} where subscription $h$ denotes the scenario index and $H$ is the total number of scenarios. 
\begin{equation}
\min_{x} \frac{1}{H}\sum_{h=1}^{H}\sum_{t=1}^{T}c_t\Delta tP_{t,h}
\label{spobj}
\end{equation}
Two types of temperature constraint formulations could be incorporated. A conservative one is formulated in \eqref{spconst1}, where the indoor air temperature of all scenarios and all steps are strictly limited within the comfort zone, referred to as SP-strict in this paper. The less conservative way is formulated in \eqref{spconst2}, where temperature violation under individual scenario is allowed, but the expected indoor air temperature overall scenarios have to satisfy the comfort zone limit. It is referred to as SP-average in this paper.
\begin{equation}
\text{SP-strict:\quad} T^{in}_{t,h} \leq T^{ub}_t, \forall t,h
\label{spconst1}
\end{equation}
\begin{equation}
\text{SP-average: \quad} \frac{1}{H}\sum_{h=1}^H T^{in}_{t,h} \leq T^{ub}_t, \forall t
\label{spconst2}
\end{equation}

The minimum up and down time constraints of SP are similar to DO's in \eqref{doformulation}, except for their dimension being expanded to all scenarios.

\subsection{Robust Optimization}\label{ro}
SP method aims to solve the robustness issue but suffers from the complexity issue. To achieve a high confidence level of feasibility and optimality, the number of scenarios required is usually large and the computational complexity of SP increases with the scenario set size. To improve the computational efficiency, the RO formulation can be used to optimize the HVAC schedule under the worst scenario of day-ahead ambient temperature within an uncertainty set. Now the objective function is modified to \eqref{roobj}, indicating a two-stage formulation. The first stage is still to minimize the total energy cost, and the second stage is to find the worst case scenario, where the power consumption reaches maximum under the worst ambient temperature $T_t^{oa*}$. Note that $P_t$ is a linear function of $T^{oa}_t$ and the first-stage decision variable $x$ is considered as a constant in the second stage, the decision variable of the second-stage problem is the ambient temperature $T^{oa}_t$. 
\begin{equation}
\min_{x} \sum_{t=1}^{T} \max_{T^{oa}_t\in \phi_t} c_t\Delta tP_{t}
\label{roobj}
\end{equation}
Note that in our paper, $\phi$ is the interval uncertainty set that varies with the predicted outdoor air temperature, i.e., $\phi_t = [\bar{\phi}_t, \underline{\phi}_t]$.

Without loss of generality, the robust counterpart (RC) of \eqref{roobj} can be reformulated to a tractable linear programming (LP) \eqref{rorc} in Appendix \ref{app1}.

\section{Distributionally Robust Optimization Formulation}\label{dro}
\subsection{DRO Formulation with a Wasserstein Ambiguity Set}\label{4a}
The SP formulation explicitly considers the uncertainty, but scalability becomes a challenge due to the number of scenarios required. In contrast, the RO formulation provides a fast and robust solution but it tends to be too conservative in practice. To overcome the shortcomings, we look at the uncertainty of entire probability distribution rather than a single prediction point using Distributionally Robust Optimization (DRO).

The proposed DRO formulation is presented below.
\begin{subequations}
\begin{align}
&\min_{x} \sum_{t=1}^{T} c_t\Delta tP_{t}\\
s.t.\quad  &\eqref{updown1}-\eqref{updown6}, \eqref{zonetemp}-\eqref{powerconsumption}\notag\\
&\max_{\mathbb{P}_t\in \mathcal{P}_t} \mathbb{E}_{\mathbb{P}_t} \left[T^{in}_t\right] \leq T^{ub}_t,\forall t \label{drocool}\\
&\mathcal{P}_t=\{\mathbb{P}_t|W\left(\mathbb{P}_t,\mathbb{Q}_t\right)\leq \varepsilon \}\label{uncertaintyset}
\end{align}
\label{droconst}
\end{subequations}

\begin{figure}[h]
{\includegraphics[width=0.37\textwidth]{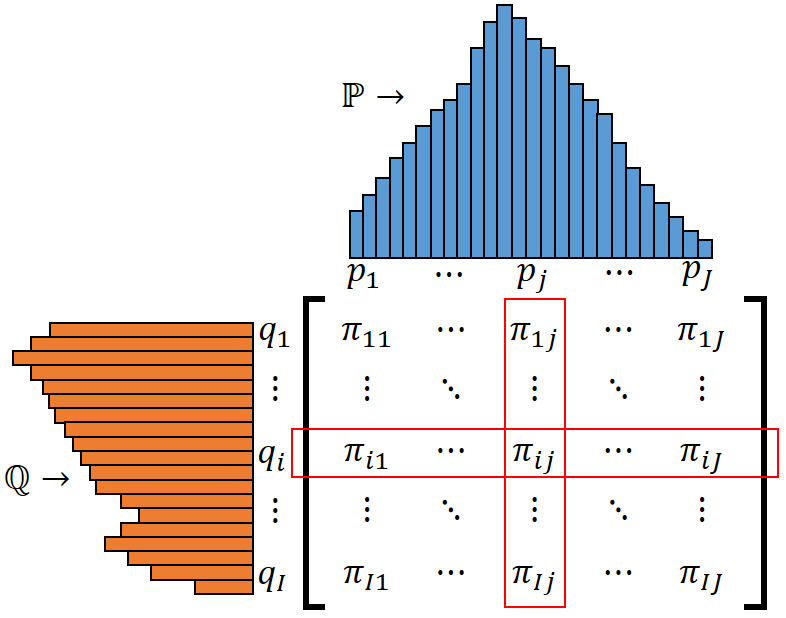}}
\caption{Wasserstein distance between two distributions and its relation to their joint distribution}
\label{wass}
\end{figure}

The initial formulation of DRO given in \eqref{droconst} is a two-stage problem, where the first-stage minimizes the total energy cost. The deterministic constraints \eqref{updown1}-\eqref{updown6}, and \eqref{zonetemp}-\eqref{powerconsumption} are implemented in the first-stage problem. The challenge lies in the second stage given by \eqref{drocool}, where the expected indoor air temperature under the worst probability distribution needs to be bounded within the comfort range. $\mathcal{P}_t$ is a Wasserstein ball (WB) ambiguity set defined in \eqref{uncertaintyset}, where $\mathbb{P}_t$ denotes an element distribution that has a Wasserstein distance to the center probability distribution $\mathbb{Q}_t$ within a radius $\varepsilon$. $W(\cdot)$ is the function of Wasserstein distance formulated in \eqref{wassformu} \cite{esfahani2018data}. 
\begin{subequations}
\begin{align}
&W\left(\mathbb{P},\mathbb{Q}\right) = \min_{\pi} \sum_{i=1}^I\sum_{j=1}^J\left| \Tilde{\xi}_{i} - \Tilde{\xi}_{j}\right|\pi_{ij}\\
&s.t.\notag\\
&\sum_{j=1}^J\pi_{ij} = q_{i},\forall i\\
&\sum_{i=1}^I\pi_{ij} = p_{j},\forall j
\end{align}
\label{wassformu}
\end{subequations}
Let $\Tilde{\xi}_i$ and $\Tilde{\xi}_j$ denote the $i^{th}$ and $j^{th}$ discretized value of $\mathbb{Q}$ and $\mathbb{P}$ and $q_i$ and $p_j$ are their probabilities respectively. Let $\pi\in \mathcal{R}^{I\times J}$ be the joint distribution matrix of the two distributions. The Wasserstein metric is the minimum effort needed to transform one distribution into the other. Note that the two probability distributions have an equal area of 1, and the difference lies in their shapes. Fig. \ref{wass} further illustrates the relationship, where the column-wise summation of $\pi$ is the probability distribution of $\mathbb{P}$, and the row-wise summation is the probability distribution of $\mathbb{Q}$. The sum of the entire $\pi$ matrix equals 1. Since $\pi$ is not unique, Wasserstein metric is to find an optimal $\pi$ that minimizes the cost of transformation.

In our case, the center distribution $\mathbb{Q}$ is the predicted empirical distribution and the radius $\varepsilon$ determines the conservativeness of the ambiguity set. Every probability distribution within the WB is a candidate and there are infinitely many of them as long as $\varepsilon > 0$. When $\varepsilon=0$, the DRO problem becomes an SP problem, otherwise the DRO problem under its initial formulation is intractable. 

To further illustrate the WB ambiguity set, let's consider an example where the empirical distribution $\mathbb{Q}_1$ is given in \eqref{q1}. $\mathbb{Q}_1$ is essentially a value prediction in the form of a probabilistic prediction. The probability of the predicted value $75$ is 100\%. 
\begin{equation}
\mathbb{Q}_1 = \left\{
\begin{aligned}
&1, T^{oa}=75^{\circ}F\\
&0, otherwise
\end{aligned}
\right.
\label{q1}
\end{equation}

Now, consider $74^{\circ}F$ and $78^{\circ}F$ as two possible temperature values in the candidate distribution $\mathbb{P}$. Given a Wasserstein distance $\varepsilon=2$, the worst probability distribution can be solved intuitively. Clearly, $78^{\circ}F$ is the value that results in a higher room temperature, therefore the higher the $p_{(78^{\circ}F)}$ is, the worse the distribution is. The global worst distribution would be $p^\dagger_{(78^{\circ}F)}=100\%$ and $p^\dagger_{(74^{\circ}F)}=0\%$. The Wasserstein distance between $\mathbb{P}^\dagger$ and $\mathbb{Q}_1$ is $(78-75)\times 100\%=3$. Such a $\mathbb{P}^\dagger$ is infeasible because it is outside of the ambiguity set. Therefore, we decrease the probability of $78^{\circ}F$, and move $\mathbb{P}$ closer to the ambiguity set. It can be concluded that the feasible worst distribution $\mathbb{P}^*$ must be sitting on the WB surface, where \eqref{intuitivesol} holds.
\begin{equation}
(78-75)\times p^*_{(78^{\circ}F)}+(75-74)\times p^*_{(74^{\circ}F)}=2
\label{intuitivesol}
\end{equation}

Considering $p^*_{(74^{\circ}F)}=1-p^*_{(78^{\circ}F)}$, we can solve for the worst probability distribution as $p^*_{(78^{\circ}F)}=50\%$ and $p^*_{(74_{\circ}F)}=50\%$. A generalized solution considering two discretized values in $\mathbb{P}$ is derived in Appendix \ref{app2}.

\subsection{Second-Stage Problem Reformulation}
Incorporating \eqref{uncertaintyset} into \eqref{drocool}, the definition of WB becomes a constraint of second-stage problem, which is given by the left-hand-side (LHS) of \eqref{wessexp1}.
\begin{equation}
  \left(
 \begin{aligned}
   &\max_{\mathbb{P}_t} \mathbb{E}_{\mathbb{P}_t} \left[T^{in}_t\right]\\
   &s.t.\\
   &W\left(\mathbb{P}_t,\mathbb{Q}_t\right)\leq \varepsilon
  \end{aligned}
  \right) \leq T^{ub}_t, \forall t
  \label{wessexp1}
\end{equation}

Substituting the room temperature $T_t^{in}$ in \eqref{wessexp1} with the ARX model \eqref{zonetemp}, the expression of $\mathbb{E}_{\mathbb{P}_t}\left[T^{in}_t\right]$ becomes \eqref{TinSamp}, where $\tilde{T}_{j,t}$ denotes the $j$th discretized value with probability $p_{j,t}$. $\tilde{T}^{in}_{j,t}$ and $\tilde{T}^{oa}_{j,t}$ are discretized values of indoor temperature and ambient air temperature probability distributions.  
\begin{subequations}
\begin{align}
\mathbb{E}_{\mathbb{P}_t} \left[T^{in}_t\right] &= \sum_{j=1}^J \Tilde{T}^{in}_{j,t}p_{j,t}\notag\\
&=\sum_{j=1}^J b_2\Tilde{T}^{oa}_{j,t}p_{j,t}+\left(b_3T^{in}_{t-1}+b_1x_t+b_0\right)\label{TinSamp}
\end{align}
\label{wessexp2}
\end{subequations}

Incorporating \eqref{wassformu} and \eqref{wessexp2} into \eqref{wessexp1}, we obtain \eqref{wessexp3}. Note the inner minimization problem is in the "$Inf\leq$" format, meaning if there exists a $\pi_{ij,t}$ that meets the condition, the constraint satisfies. Therefore, minimization can be removed and expression \eqref{wessexp3} can be simplified to regular constraints in \eqref{wessexp} \cite{esfahani2018data}.
\begin{equation}
  \left(
 \begin{aligned}
   &\max_{\mathbb{P}_t}\sum_{j=1}^J b_2\Tilde{T}^{oa}_{j,t}p_{j,t}+\left(b_3T^{in}_{t-1}+b_1x_t+b_0\right)\\
   &s.t.\\
   &\left(
   \begin{aligned}
     &\min_{\pi_{ij,t}} \sum_{i=1}^I\sum_{j=1}^J\left| \Tilde{T}^{oa}_{i,t} - \Tilde{T}^{oa}_{j,t}\right|\pi_{ij,t}\\
     &s.t.\\
     &\sum_{j=1}^J\pi_{ij,t} = q_{i,t},\forall i\\
     &\sum_{i=1}^I\pi_{ij,t} = p_{j,t},\forall j\\
     \end{aligned}\right)\leq \varepsilon
  \end{aligned}
  \right) \leq T^{ub}_t, \forall t
  \label{wessexp3}
\end{equation}
\begin{equation}
  \left(
 \begin{aligned}
   &\max_{\mathbb{P}_t}\sum_{j=1}^J b_2\Tilde{T}^{oa}_{j,t}p_{j,t}+\left(b_3T^{in}_{t-1}+b_1x_t+b_0\right)\\
   &s.t.\\
   &\sum_{i=1}^I\sum_{j=1}^J\left| \Tilde{T}^{oa}_{i,t} - \Tilde{T}^{oa}_{j,t}\right|\pi_{ij,t}\leq \varepsilon\\
   &\sum_{j=1}^J\pi_{ij,t} = q_{i,t},\forall i\\
   &\sum_{i=1}^I\pi_{ij,t} = p_{j,t},\forall j\\
  \end{aligned}
  \right) \leq T^{ub}_t, \forall t
  \label{wessexp}
\end{equation}

By now, the proposed second-stage DRO constraint \eqref{drocool} is transformed into the linear constraint in \eqref{wessexp}, where the goal is to find the worst probability distribution $\mathbb{P}_t$ that is bounded within the WB of the empirical probability distribution $\mathbb{Q}_t$. 
$T_{t-1}^{in}$ and $x_t$ are not second-stage decision variables, therefore can be regarded as constants.

\subsection{Tractable Reformulation using Dualization}
The LHS of \eqref{wessexp} can be further reformulated into \eqref{wessexp4} by substituting $p_{j,t}$ in the objective function with $\sum_{i=1}^I\pi_{ij,t}$ to cancel out $p_{j,t}$. Then $\pi_{ij,t}$ becomes the only decision variable of the second-stage problem.
\begin{subequations}
 \begin{align}
   &\max_{\pi_t}\sum_{j=1}^J\sum_{i=1}^I b_2\Tilde{T}^{oa}_{j,t}\pi_{ij,t} + \left(b_3T^{in}_{t-1}+b_1x_t+b_0\right)\\
   &s.t.\notag\\
   &\sum_{i=1}^I\sum_{j=1}^J\left| \Tilde{T}^{oa}_{i,t} - \Tilde{T}^{oa}_{j,t}\right|\pi_{ij,t}\leq \varepsilon \quad \left(\lambda_t\right)\label{lambda}\\
   &\sum_{j=1}^J\pi_{ij,t} = q_{i,t},\forall i \quad \left(s_{i,t}\right)\label{si}
  \end{align}
\label{wessexp4}
\end{subequations}
Let $\lambda_t$ and $s_{i,t}$ be the Lagrangian multipliers, the Lagrangian function of \eqref{wessexp4} is the $\mathcal{L}_t\left(\pi_t,\lambda_t,s_t\right)$ formulated in \eqref{lagrangian}, which can be further reformulated into \eqref{lagrangian2} by aggregating the $\pi_{ij,t}$ terms. 
\begin{equation}
\begin{aligned}
\mathcal{L}_t\left(\pi_t,\lambda_t,s_t\right)=&\sum_{j=1}^J\sum_{i=1}^I b_2\Tilde{T}^{oa}_{j,t}\pi_{ij,t} + \left(b_3T^{in}_{t-1}+b_1x_t+b_0\right)\\
&+ \left(\varepsilon - \sum_{i=1}^I\sum_{j=1}^J\left| \Tilde{T}^{oa}_{i,t} - \Tilde{T}^{oa}_{j,t}\right|\pi_{ij,t}\right)\lambda_t \\
&+ \sum_{i=1}^I\left(q_{i,t} - \sum_{j=1}^J\pi_{ij,t}\right)s_{i,t} 
\end{aligned}
\label{lagrangian}
\end{equation}
\begin{equation}
\begin{aligned}
\mathcal{L}_t\left(\pi_t,\lambda_t,s_t\right)=&\varepsilon\lambda_t + \sum_{i=1}^Iq_{i,t}s_{i,t} + \left(b_3T^{in}_{t-1}+b_1x_t+b_0\right)\\
&+ \sum_{j=1}^J\sum_{i=1}^I \left(b_2\Tilde{T}^{oa}_{j,t} - \lambda_t\left| \Tilde{T}^{oa}_{i,t} - \Tilde{T}^{oa}_{j,t}\right| - s_{i,t}\right)\pi_{ij,t} 
\end{aligned}
\label{lagrangian2}
\end{equation}
Considering $\pi_{ij,t}\geq0$, there exists a bounded supremum only if the coefficient of $\pi_{ij,t}$ is non-positive, i.e. $\left| \Tilde{T}^{oa}_{i,t} - \Tilde{T}^{oa}_{j,t}\right|\lambda_t + s_{i,t}\geq b_2\Tilde{T}^{oa}_{j,t}$. The supremum $p^* = \varepsilon\lambda_t + \sum_{i=1}^Iq_{i,t}s_{i,t} + \left(b_3T^{in}_{t-1}+b_1x_t+b_0\right)$ is achieved at $\pi_{ij,t}=0$. According to the strong duality theorem \cite{boyd2004convex}, the Lagrangian function becomes a function of $\lambda_t$ and $s_{i,t}$, and its infimum is equal to the supremum of the original problem. Hence, the maximization problem in \eqref{wessexp4} is dualized into a minimization problem \eqref{wessexp5}. Together with the upper bound of the temperature setpoint, constraint \eqref{wessexp} becomes the inequality \eqref{wessexp5}. 
\begin{equation}
  \left(
 \begin{aligned}
   &\min_{\lambda_t,s_t}\varepsilon\lambda_t + \sum_{i=1}^Iq_{i,t}s_{i,t} + \left(b_3T^{in}_{t-1}+b_1x_t+b_0\right)\\
   &s.t.\\
   &\left| \Tilde{T}^{oa}_{i,t} - \Tilde{T}^{oa}_{j,t}\right|\lambda_t + s_{i,t}\geq b_2\Tilde{T}^{oa}_{j,t}, \forall i,j\\
  \end{aligned}
  \right) \leq T^{ub}_t, \forall t
  \label{wessexp5}
\end{equation}
Similar to the sub-problem in \eqref{wessexp3}, the "min" can be removed due to "$Inf \leq$" format. Hence, the final tractable formulation of the proposed DRO under cooling conditions is given by \eqref{coolfinal}.
\begin{subequations}
\begin{align}
&\min_{x} \sum_{t=1}^{T} c_t\Delta tP_{t}\\
&s.t.\notag\\
&\varepsilon\lambda_t + \sum_{i=1}^Iq_{i,t}s_{i,t} + \left(b_3T^{in}_{t-1}+b_1x_t+b_0\right)\leq T^{ub}_t, \forall t\label{wessexp6}\\
&\left| \Tilde{T}^{oa}_{i,t} - \Tilde{T}^{oa}_{j,t}\right|\lambda_t + s_{i,t}\geq b_2\Tilde{T}^{oa}_{j,t}, \forall i,j,t\label{wessexp7}\\
&\eqref{updown1}-\eqref{updown6}, \eqref{zonetemp}-\eqref{powerconsumption}\notag
\end{align}
\label{coolfinal}
\end{subequations}

\section{Case Study}\label{case}
Two test cases are carried out in this section. The first case is a simple test case that provides an intuitive explanation of the proposed DRO formulations. Then the DRO approach is applied to optimize a one-day schedule of a commercial building, with the prediction of the ambient temperature following Gaussian distribution. The optimization solver used in this section is Gurobi, and the environment is Python.

\subsection{An Intuitive Example}
This subsection aims to provide a straightforward test case following the example in Section \ref{4a} with $\mathbb{Q}_1$ being the ambient temperature forecast of $75^{\circ}F$. The example given is a single-step decision-making considering two possible ambient temperature values. The power consumption model and room temperature model incorporated in this case are formulated in \eqref{simplecasemodels}. The room temperature at the previous step is $T^{in}_0=76$, and the upper bound of comfort zone is set to $76^{\circ}F$. Then the initial formulation of this case can be formulated as \eqref{simpleformu}, where the $WB(\mathbb{Q}_1,2)$ denotes the Wasserstein ball with the center of $\mathbb{Q}_1$ and a radius of $\varepsilon=2$.
\begin{subequations}
\begin{align}
&P=100x + 0.3T^{oa}\\
&T^{in}=-3x + 0.3T^{oa} + 0.7T^{in}_0
\end{align}
\label{simplecasemodels}
\end{subequations}
\begin{subequations}
\begin{align}
&\min_{x} 0.1\left(100x + 0.3T^{oa}\right)\\
&s.t.\notag\\
&\max_{\mathbb{P}\in WB\left(\mathbb{Q}_1,2\right)} \mathbb{E}_{\mathbb{P}} \left[-3x + 0.3T^{oa} + 0.7T^{in}_0\right] \leq 76\label{20b}
\end{align}
\label{simpleformu}
\end{subequations}

The results of 7 scenarios with prediction of $\mathbb{Q}_1$ are summarized in Table \ref{onesteptwotempq1}, where $\mathbb{E}\left[T^{in}\mid_{x=1}\right]$ and $\mathbb{E}\left[T^{in}\mid_{x=0}\right]$ are the expected value of room temperature with and without HVAC under the worst case distribution which is shown in columns 2 and 3. $x^*$ is the optimal HVAC on/off status solved by the tractable formulation of the proposed DRO method. Take the first case for example. When two discretized value $75^{\circ}F$ and $77^{\circ}F$ are considered, the worst probability distribution is when $p_{77^{\circ}F}=1$. Inserting this probability into \eqref{20b}, we find that the expected indoor temperature is $76.3^{\circ}F$ when HVAC is off and $73.3^{\circ}F$ when HVAC is on. To avoid temperature violation, HVAC needs to be on $(x=1)$, which is exactly the solution from tractable reformulation \eqref{coolfinal}. The solutions from all seven cases can be manually computed to confirm the correctness of tractable reformulation \eqref{coolfinal} from the very original two-stage DRO formulation \eqref{droconst}.  It is also noted that different discretization values impact the optimal solutions. 
\begin{table}[ht!]
\begin{center}
\caption{Intuitive example result with $\mathbb{Q}_1$}
\label{onesteptwotempq1}
\scalebox{0.9}{
\begin{tabular}{ c c c c c c}
\hline
\textbf{Index} & \textbf{$\Tilde{\xi}_1(p_{\Tilde{\xi}_1})$}& \textbf{$\Tilde{\xi}_2(p_{\Tilde{\xi}_2})$} & \textbf{$\mathbb{E}\left[T^{in}\mid_{x=0}\right]$} & \textbf{$\mathbb{E}\left[T^{in}\mid_{x=1}\right]$}&\textbf{$x^*$}\\
\hline
1 & 75(0\%) & 77(100\%) & 76.3& 73.3 & 1 \\
2 & 74(50\%) & 78(50\%) & 76 & 73 & 0 \\
3 & 75(33.3\%) & 78(66.70\%) & 76.3 & 73.3 & 1 \\
4 & 76(50\%) & 78(50\%) & 76.3 & 73.3 & 1 \\
5 & 74(66.7\%) & 79(33.3\%) & 75.9 & 72.9 & 0 \\
6 & 75(50\%) & 79(50\%) & 76.3 & 73.3 & 1 \\
7 & 76(66.70\%) & 79(33.3\%) & 76.3 & 73.3 & 1 \\
\hline
\end{tabular}}
\end{center}
\end{table}

\subsection{Practical Test Case}

The practical test case utilizes real measurement data collected from a commercial building to learn the zone temperature and power consumption model. Table \ref{case2para} summarizes the coefficients of the models used in this case. 
\begin{table}[ht!]
\begin{center}
\caption{Parameters of the practical test case}
\label{case2para}
\begin{tabular}{ c | c || c | c}
\hline
\textbf{Parameter} & \textbf{Value} & \textbf{Parameter} & \textbf{Value}\\
\hline
$b_1$& -2.07 & $a_1$& 70.7\\
$b_2$& 0.15 & $a_2$& 0.24\\
$b_3$& 0.45 & $a_0$& -17.8\\
$b_0$& 37.9 & $\Delta t$& 0.1\\
$T^{in}_0$& $80^{\circ}F$ &  &  \\
\hline
\end{tabular}
\end{center}
\end{table}
The DO, SP, RO, and the proposed DRO methods are implemented to solve the one-day schedule of the building. The obtained optimal schedules are tested under two testing sets, the regular scenario testing set and the extreme scenario testing set, to validate their optimality and robustness. Each testing set contains 1000 randomly generated scenarios. The scenarios in the regular testing set are sampled from the predicted Gaussian distribution, while the scenarios in the extreme testing set are sampled from a set of distributions randomly, including Gaussian distribution with inaccurate mean and standard deviation, uniform distribution with random ranges, and beta distribution with random parameters of $\alpha$ and $\beta$, to mimic the unknown nature of underlying probability distributions of ambient temperature.

\begin{figure}[h]
{\includegraphics[width=0.5\textwidth]{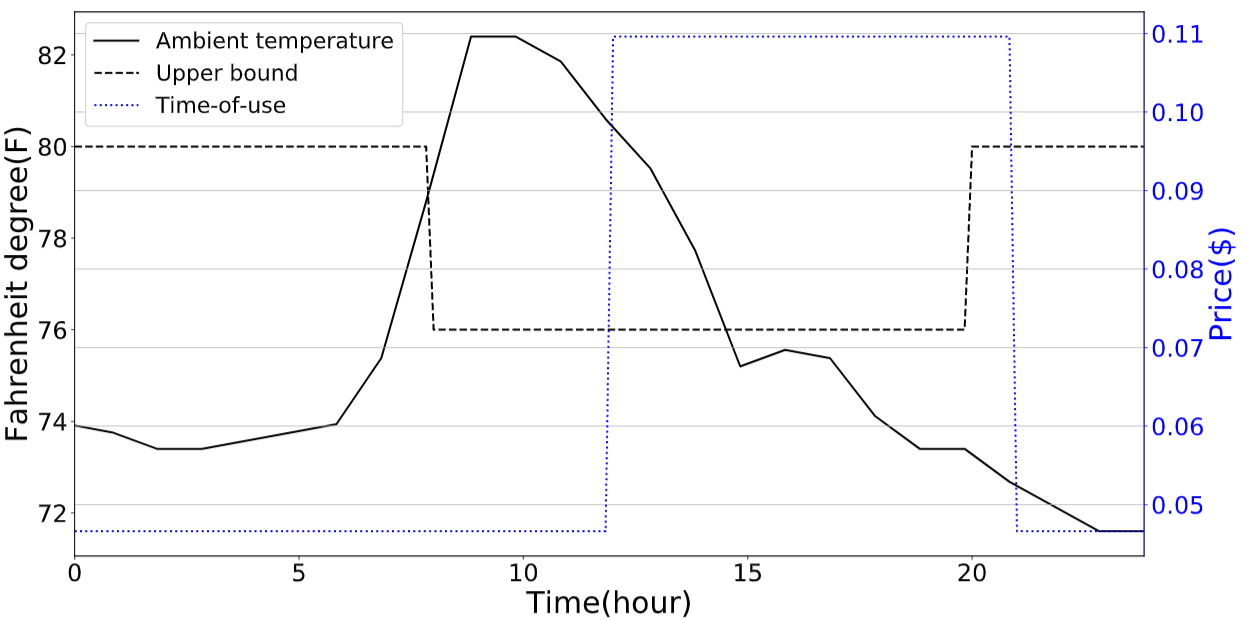}}
\caption{Trajectories of predicted outdoor air temperature, zone temperature upper bound, and TOU price}
\label{case1settings}
\end{figure}

The optimality evaluation metric is the total cost. The robustness evaluation metrics are the number of steps with temperature violations $\left(V^{num}\right)$, and the temperature violation mileage $\left(V^{mil}\right)$, defined in \eqref{vio}, where $\mathbbm{1}(\cdot)$ is the indicator function.
\begin{subequations}
\begin{align}
&V^{num}=\sum_{t=1}^T \mathbbm{1}_{\left(T^{in}_t-T^{ub}_t\right)}\\
&V^{mil}=\sum_{t=1}^T \mathbbm{1}_{\left(T^{in}_t-T^{ub}_t\right)}\times \left(T^{in}_t-T^{ub}_t\right)    
\end{align}
\label{vio}
\end{subequations}
The predicted outdoor air temperature is assumed to follow Gaussian distributed in the day-ahead prediction, and the predicted mean value $\mu_t$ over the day is indicated by the solid line in Fig. \ref{case1settings}. The predicted standard deviation of $T^{oa}$ is set to $0.5^{\circ}F$ for all steps. The probabilistic prediction of outdoor air temperature at step $t$ is $\mathbb{Q}_t= N(\mu_t,0.5^2)$. The dashed line in Fig. \ref{case1settings} denotes the zone temperature upper bound of the day, which is $76^{\circ}F$ for the working hours (8 am to 8 pm) and $80^{\circ}F$ for the off-work hours (12 am to 8 am and 8 pm to 12 am). The dotted line denotes the aggregated TOU price, including the on-peak (12 pm to 9 pm) and off-peak (12 am to 12 pm and 9 pm to 12 am) prices of fuel cost, electricity cost, and other miscellaneous costs charged by utility companies. \cite{duke}.

\begin{figure}[h]
{\includegraphics[width=0.5\textwidth]{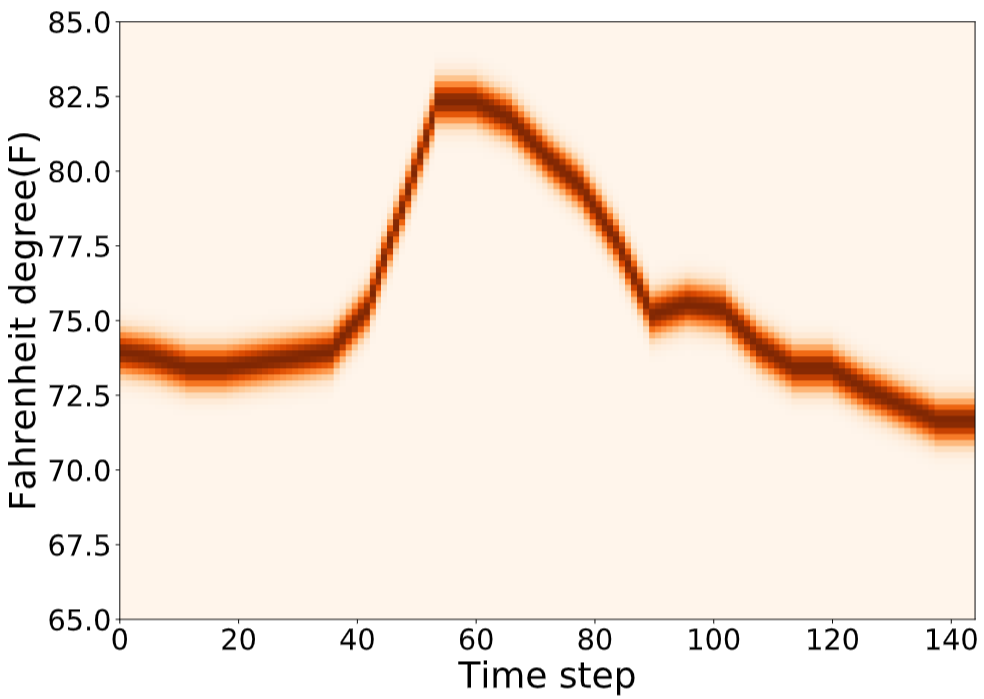}}
\caption{Probabilistic prediction of outdoor air temperature}
\label{case1probpred}
\end{figure}

Fig. \ref{case1probpred} shows the probabilistic prediction of the entire day, where darker color indicates higher probability. The time resolution is 10 minutes, so the total number of step $T$ in a day is $24\times 6=144$. 

The feasible region of $T^{oa}$ is set to $[65,85]$, which has at least $6 \sigma$-wide margins on both sides to the range of the predicted mean, therefore can be considered sufficient for discretization. The feasible region is equally discretized into 100 segments to obtain the discretized values $\tilde{T}_t^{oa}$. 

\begin{figure}[h]
{\includegraphics[width=0.5\textwidth]{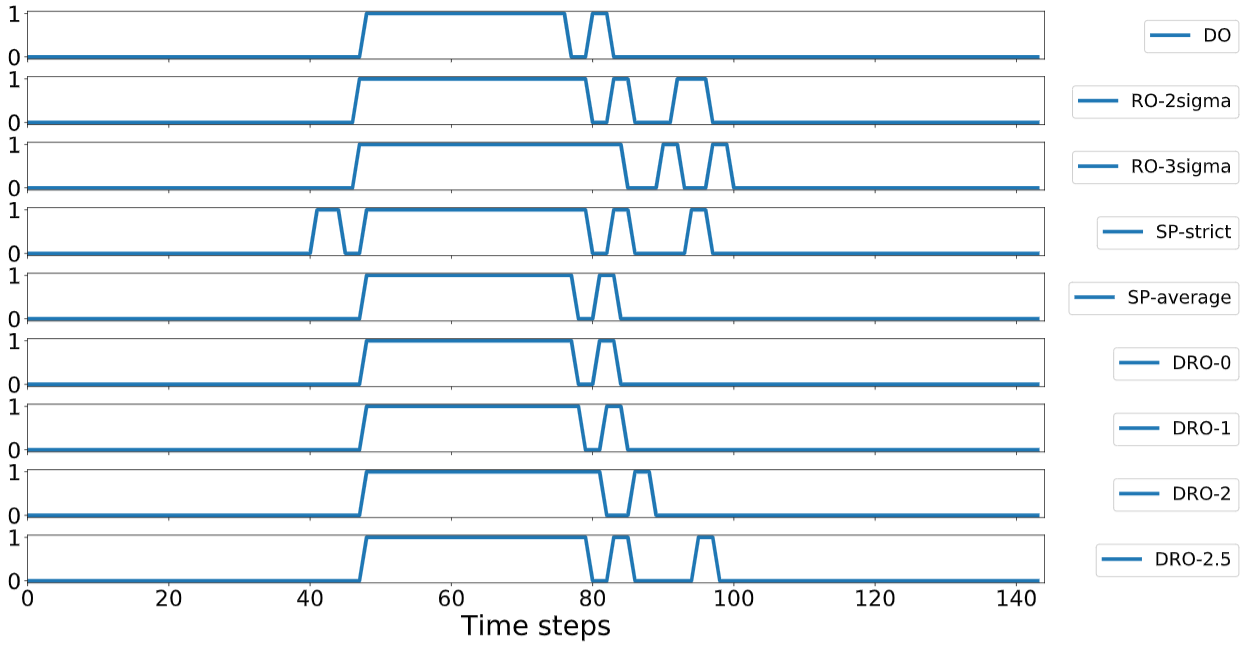}}
\caption{Optimized schedules}
\label{case1optschedules}
\end{figure}

Fig. \ref{case1optschedules} shows the optimized schedules, where DO is solved according to the predicted mean, RO-2$\sigma$ is solved with 95\% confidence interval, of which the robust feasible region at each step is $\phi_t=[\mu_t-2\sigma,\mu_t+2\sigma]$. Similarly, RO-3$\sigma$ is solved with 99.7\% confidence interval, of which the robust feasible region at each step is $\phi_t=[\mu_t-3\sigma,\mu_t+3\sigma]$. SP-strict and SP-average are the solutions of SP formulations with temperature constraints defined upon individual scenarios and the average of all scenarios. DRO-0, DRO-1, DRO-2, and DRO-2.5 are the solutions of the proposed DRO methods with Wasserstein radius $\varepsilon$ equal to 0, 1, 2 and 2.5 respectively. It can be observed that the discrepancies among these schedules occur in the early morning around 6-8 am and afternoon around 12-4 pm. The former time span involves the ambient temperature ramping up and room temperature upper bound switching to $76^{\circ}F$, and the later time span involves the TOU step-up and ambient temperature ramping down. 

\begin{table}[h]
\begin{center}
\caption{Testing results of the optimized schedules}
\label{case1comp}
\scalebox{1}{
\begin{tabular}{ c c | c c c c c}
\hline
\textbf{Index} & \textbf{Method} & \textbf{Cost(\$)}& \textbf{$V^{num}_{reg}$}& \textbf{$V^{mil}_{reg}$}& \textbf{$V^{num}_{ext}$}& \textbf{$V^{mil}_{ext}$}\\
\hline
1&DO & 71.506& 0.771& 0.125& 12.27& 3.312\\
2&RO-2$\sigma$ & 86.896& 0& 0& 0.29& 0.019\\
3&RO-3$\sigma$ & 92.376& 0& 0& 0& 0\\
4&SP-strict & 86.349& 0& 0& 0.37& 0.032\\
5&SP-average & 73.333& 0.201& 0.021& 9.66& 2.292\\
6&DRO-0 & 73.333& 0.201& 0.021& 9.66& 2.292\\
7&DRO-1 & 75.159& 0.027& 0.003& 7.5& 1.754\\
8&DRO-2 & 80.639& 0.003& 0.0002& 4.92& 1.168\\
9&DRO-2.5 & 82.466& 0& 0& 0.32& 0.024\\
\hline
\end{tabular}}
\end{center}
\end{table}

Table \ref{case1comp} summarizes the performance of the 9 optimal schedules, where the subscripts $reg$ and $ext$ denotes the regular testing set and extreme testing set respectively. It can be observed that DO has the lowest operating cost, but comes with high violations in both testing sets. Its poor robustness is because DO only takes the predicted mean into consideration, but in reality, predictions are never perfect. The two RO solutions with different confidence intervals show strong robustness, especially RO-3$\sigma$ has no violations observed in both testing sets, but it also has the highest cost due to its conservativeness. The robustness of SP-strict is close to RO-2$\sigma$ in that their violation number and mileage are at the same level in both testing sets. Their costs are approximately the same, and the cost of SP-strict is 0.6\% lower than RO-3$\sigma$. 

The solution of SP-average is less conservative than SP-strict because its temperature constraint is based on the average of all scenarios. The schedule of DRO-0 and SP-average are the same, and they have the same cost and robustness. The DRO method with an ambiguity set radius 0, is essentially an SP problem because the ambiguity set only includes the predicted empirical distribution. The only difference is SP-average is solved in a scenario-based way and DRO-0 is solved in the discretized LP manner. From DRO-0 to DRO-2.5, the conservativeness and cost increase with the increase of $\varepsilon$, as a result, the violation number and mileage decrease. From the feasibility perspective, DRO-2 is already very close to SP-strict and RO-2$\sigma$ in the regular test set, and DRO-2.5 yields similar robustness to SP-strict and RO-2$\sigma$ in the extreme test set. The violation number and mileage of DRO-2.5 even become slightly lower than SP-strict. In terms of optimality, the cost of DRO-2 is 6.6\% lower than SP-strict and RO-2$\sigma$. DRO-2.5 is more expensive than DRO-2 by 2.3\% but is cheaper than SP-strict and RO-2$\sigma$ by 4.5\%. The cost of DRO solutions with $\varepsilon\in\left(0,2.5\right]$ is bounded by the cost of DRO-0 and DRO-2.5. It can be concluded that, with the proper selection of $\varepsilon$, the proposed DRO method can provide competitively robust schedules compared to SP and RO methods with lower costs. The proposed DRO method provides building managers with granular options according to their trade-off between economic and reliability preferences.

\section{Conclusions}\label{conclusions}
This paper studies the formulations of HVAC setpoint optimization methods, including DO, SP, and RO. A DRO approach based on the Wasserstein ambiguity set is then proposed to enhance the robustness of the optimized schedules by considering the uncertainty of ambient temperature predictions. The DRO method minimizes the total operation cost, while the room temperature under the worst probability distribution still satisfies occupant comfort range. The DRO problem is reformulated into a tractable MILP form and implemented in two test cases. The numerical results indicate that DRO can yield comparable robustness with SP and RO, but with much lower cost. The proposed DRO approach provides granular options regarding the risk-benefit preference. The future work is to explore the uncertainties caused by both input data and building models.

\begin{appendices}
\section{Robust Counterpart (RC)}\label{app1}
The RO formulation \eqref{roobj} is equivalent to \eqref{ro2} by introducing an auxiliary variable $z_t$ that denotes the upper bound of energy cost.
\begin{subequations}
\begin{align}
&\min_{x,z} \sum_{t=1}^{T} z_t\\
s.t.\quad &\max_{\substack{T^{oa}_{t}\in \phi_t\\x_t}} c_t\Delta tP_{t}\leq z_t, \forall t 
\end{align}
\label{ro2}
\end{subequations}

Substituting $P_t$ with the HVAC power consumption model \eqref{powerconsumption}, the Robust Counterpart (RC) becomes \eqref{ro3}, where the second-stage problem is a linear function of $T^{oa}_t$ and its feasible region is the interval uncertainty set $\phi_t$.  
\begin{subequations}
\begin{align}
&\min_{x,z} \sum_{t=1}^{T} z_t\\
&s.t.\notag \\
&\max_{\substack{T^{oa}_{t}\in \phi_t\\x_t}} c_t\Delta t\left(a_{1}x_{t} + a_{2}T^{oa}_t + a_0\right)\leq z_t, \forall t 
\end{align}
\label{ro3}
\end{subequations}

Consider that $x$ is a binary variable, the second-stage problem is a bi-linear problem and its supremum must be achieved at the boundary conditions. The tractable formulation for RC is given in \eqref{rorc}, where \eqref{ro4} and \eqref{ro5} denotes the boundary conditions when $x=1$ and \eqref{ro6} and \eqref{ro7} denotes that of $x=0$.
\begin{subequations}
\begin{align}
&\min_{z} \sum_{t=1}^{T} z_t\\ 
&s.t. \notag\\
&z_t \geq c_t\Delta t \left( a_{2}\bar{\phi}_t + a_0\right) + c_t\Delta ta_{1} , \forall t\label{ro4}\\
&z_t \geq c_t\Delta t \left( a_{2}\underline{\phi}_t + a_0\right) + c_t\Delta ta_{1} , \forall t\label{ro5}\\
&z_t\geq c_t\Delta t \left( a_{2}\bar{\phi}_t + a_0\right), \forall t\label{ro6}\\
&z_t\geq c_t\Delta t \left( a_{2}\underline{\phi}_t + a_0\right), \forall t\label{ro7}
\end{align}
\label{rorc}
\end{subequations}

\section{Worst Case Distribution}\label{app2}
According to the definition of Wasserstein metric-based ambiguity set in \eqref{uncertaintyset} and \eqref{wassformu}, the worst case distribution of two possible value satisfies \eqref{worstcasep1}, where $\Tilde{\xi}_1\leq\Tilde{\xi}_2$. $p_{\Tilde{\xi}_1}$ and $p_{\Tilde{\xi}_2}$ are their associated probabilities. 
\begin{equation}
\left|\Tilde{\xi}_1-75\right|p_{\Tilde{\xi}_1}+\left|\Tilde{\xi}_2-75\right|p_{\Tilde{\xi}_2}\leq 2
\label{worstcasep1}
\end{equation}

Considering $p_{\Tilde{\xi}_1} + p_{\Tilde{\xi}_2} =1$, we can get \eqref{worstcasep2} by substituting $p_{\Tilde{\xi}_1}=1-p_{\Tilde{\xi}_2}$ into \eqref{worstcasep1}.
\begin{equation}
\left(\left|\Tilde{\xi}_2-75\right|-\left|\Tilde{\xi}_1-75\right|\right)p_{\Tilde{\xi}_2}\leq 2-\left|\Tilde{\xi}_1-75\right|
\label{worstcasep2}
\end{equation}

The feasible region of $p_{\Tilde{\xi}_2}$ is different under the following three scenarios:
\begin{enumerate}
    \item When $\Tilde{\xi}_2\geq \Tilde{\xi}_1\geq 75$, we have $\left|\Tilde{\xi}_1-75\right|\geq0$, $\left|\Tilde{\xi}_2-75\right|\geq0$, and $\left|\Tilde{\xi}_2-75\right|-\left|\Tilde{\xi}_1-75\right|\geq0$. The solution of \eqref{worstcasep2} is $p_{\Tilde{\xi}_2}\leq \frac{77-\Tilde{\xi}_1}{\Tilde{\xi}_2-\Tilde{\xi}_1} \text{ and } p_{\Tilde{\xi}_2}\leq1$.
    \item When $\Tilde{\xi}_2\geq 75 \geq \Tilde{\xi}_1$, and $\Tilde{\xi}_1+\Tilde{\xi}_2\geq 150$, we have $\left|\Tilde{\xi}_1-75\right|\geq0$, $\left|\Tilde{\xi}_2-75\right|\leq0$, and $\left|\Tilde{\xi}_2-75\right|-\left|\Tilde{\xi}_1-75\right|\geq0$. The solution of \eqref{worstcasep2} is $p_{\Tilde{\xi}_2}\leq \frac{\Tilde{\xi}_1-73}{\Tilde{\xi}_2-\Tilde{\xi}_1} \text{ and } p_{\Tilde{\xi}_2}\leq1$.
    \item When $\Tilde{\xi}_2\geq 75 \geq \Tilde{\xi}_1$, and $\Tilde{\xi}_1+\Tilde{\xi}_2\leq 150$, we have $\left|\Tilde{\xi}_1-75\right|\geq0$, $\left|\Tilde{\xi}_2-75\right|\leq0$, and $\left|\Tilde{\xi}_2-75\right|-\left|\Tilde{\xi}_1-75\right|\leq0$. The solution of \eqref{worstcasep2} is $1\geq p_{\Tilde{\xi}_2}\geq \frac{\Tilde{\xi}_1-73}{\Tilde{\xi}_2-\Tilde{\xi}_1}$.
\end{enumerate}

It can be seen from \eqref{TinSamp} that, when the discretized value is fixed, the expected room temperature $\mathbb{E}\left[T^{in}\right]$ is linear to $p$. The worst case distribution  that maximize $\mathbb{E}\left[T^{in}\right]$ must be achieved at the upper bound of $p_{\Tilde{\xi}_2}$ and the lower bound of $p_{\Tilde{\xi}_1}$, because the former one has a larger positive coefficient. Hence, the solution of the the worst distribution $\mathbb{P}^*$ under two possible values can be summarized as \eqref{worstcasep}. 
\begin{equation}
\mathbb{P}^*=\left\{
\begin{aligned}
&p_{\Tilde{\xi}_1}=1-p_{\Tilde{\xi}_2}\\
&p_{\Tilde{\xi}_2}=
\left\{
\begin{aligned}
&\min\left(\frac{77-\Tilde{\xi}_1}{\Tilde{\xi}_2-\Tilde{\xi}_1},1\right), \Tilde{\xi}_2\geq \Tilde{\xi}_1\geq 75\\
&\min\left(\frac{\Tilde{\xi}_1-73}{\Tilde{\xi}_2-\Tilde{\xi}_1},1\right), \begin{aligned}&\Tilde{\xi}_2\geq 75 \geq \Tilde{\xi}_1\\
&\text{and}\\
&\Tilde{\xi}_1+\Tilde{\xi}_2\geq 150\end{aligned}\\
&1,\Tilde{\xi}_2\geq 75 \geq \Tilde{\xi}_1 \text{ and } \Tilde{\xi}_1+\Tilde{\xi}_2\leq 150
\end{aligned}
\right.
\end{aligned}
\right.
\label{worstcasep}
\end{equation}

Note that all feasible values are covered by the three scenarios in \eqref{worstcasep}, because $\Tilde{\xi}_1$ and $\Tilde{\xi}_2$ must satisfy $\Tilde{\xi}_2\geq76$ and $\min\left(\left|\Tilde{\xi}_1-75\right|,\left|\Tilde{\xi}_2-75\right|\right)\leq2$ to ensure meaningful temperature constraint and non-empty ambiguity set. 

\end{appendices}

\bibliographystyle{IEEEtran}
\bibliography{DRO.bib}

\end{document}